\documentclass[12pt]{article}

\usepackage[margin=3cm]{geometry}
\usepackage{amsfonts,amssymb,amsmath,epsfig,euscript}
\usepackage{hyperref}
\usepackage{breakurl}
\usepackage{fancyhdr}
\usepackage{setspace}
\usepackage{graphicx}
\usepackage{CJKutf8} 

\newcommand{\arxiv}[1]{\href{http://arxiv.org/abs/#1}{\tt arXiv:\nolinkurl{#1}}}

\title{Introducing the Sangaku Archive}
\author{David Clark\thanks{Research for this work was partially conducted while a Visiting Scholar at the Waseda Institute for Advanced Study, with support from the Rashkind Family Endowment.}}

\begin{document}

\maketitle

The 1641 edition of Yoshida Mitsuyoshi's \emph{Jink\={o}ki}
\begin{CJK}{UTF8}{min}
(塵劫記),
\end{CJK}
a math textbook that had already enjoyed more than a decade of feverish popularity, concluded with something new: a list of twelve unsolved geometry problems. Intended as challenges for readers, these problems, which came to be known as \emph{idai}
\begin{CJK}{UTF8}{min}
(遺題
\end{CJK}
or ``bequeathed problems"), began as riffs on classical exercises imported from China. The practice of publishing idai quickly caught fire among math enthusiasts, where a culture of competition developed alongside a desire to push beyond Chinese hand-me-downs. The great mathematician Seki Takakazu (?--1708) cut his teeth on a collection of bequeathed problems from 1671, though he later criticized the phenomenon for its tendency toward haste and one-upmanship \cite{horiuchi1}. 

It was within this milieu that the tradition of posting \emph{sangaku} tablets arose. Sangaku
\begin{CJK}{UTF8}{min}
(算額, 
\end{CJK}
literally ``framed calculations")---wooden tablets inscribed with geometry problems and posted at Buddhist temples and Shinto shrines---were an offshoot of an older tradition in which patrons would present votive tablets (\emph{ema} or
\begin{CJK}{UTF8}{min}
絵馬)
\end{CJK}
as religious offerings. Sangaku problems were not all strictly idai, since most claimed to be solved; the tablets would state a numerical answer (if one was desired) as well as a formula or technique used, but not a complete solution. (Much has been written about sangaku; see \cite{fukrot1,fukrot2, clark1,hosking1}.)

Sangaku were originally products of the Edo period (1603--1868), a time of international isolation and domestic peace during which art forms like flower arranging, haiku poetry, kabuki theater, and \emph{ukiyo-e} woodblock printing flourished. Traditional Japanese mathematics (also known today as \emph{wasan} or
\begin{CJK}{UTF8}{min}
和算)
\end{CJK}
developed within this cultural renaissance, and one can view sangaku tablets---meticulously laid-out, ornately framed, and frequently brightly colored---as artistic outputs of the discipline. This understates their purpose, however. Competition was fierce among Edo period mathematics schools, which frequently used sangaku as advertisements for prospective pupils \cite{horiuchi2}. Implicitly, they would say: ``Look at these problems our students solved; can you solve them?"

Thus, sangaku were media devices serving a variety of roles: religious offerings, works of art, billboards, and vectors for mathematical communication. Their breadth of purpose was matched by their geographic distribution: many thousands of sangaku tablets were dedicated throughout the Japanese archipelago, at temples and shrines in big cities and rural villages. Made of wood and frequently installed outside, most original tablets have been lost, though the contents of many were recorded in travel journals and other documents.

There are roughly 900 surviving sangaku tablets, scattered across Japan; many are still in the possession of the temple or shrine in which they were first hung. Until recently, the only digital repository of this artifact collection was hosted not by a university library or national museum, but rather by retired teacher and amateur scholar Hiroshi Kotera. Kotera's primarily Japanese-language website, \url{www.wasan.jp}, contains an enormous amount of valuable wasan-related information, but the sangaku repository itself is difficult to navigate: created using static HTML, it is not searchable, and its photographs have uneven quality and often load slowly.

Cultural preservationists, while honoring Kotera's diligent and painstaking work, felt it was time for an upgrade. Thus began the Sangaku Archive project.

\bigskip
\centerline{\rule{2in}{0.4pt}}
\bigskip

Intellectual historian Antonia Karaisl initiated the project in 2023 as a researcher at Waseda University in Tokyo. Her objective was to create a bilingual, searchable database of extant sangaku tablets---in addition to written records of lost tablets---with sharp, scalable images. Employing the latest tools of the digital humanities, the Sangaku Archive would be built using the International Image Interoperability Framework (IIIF), which efficiently standardizes delivery of information across browsers. In IIIF, each image is paired with a customizable list of searchable metadata; for sangaku tablets, these metadata include the year and place of dedication, the name of the school and/or individual(s) who dedicated the tablet, and the tablet's physical dimensions. IIIF also supports deep zooming of images, allowing users to fully and efficiently explore visual objects without the need for heavy computation or a fast internet connection \cite{IIIF}.

Deep zooming requires high-quality photographs; in 2023, these did not exist for most sangaku tablets. This would be a major task in building the archive. There are extant sangaku in more than three-quarters of Japan's 47 prefectures. In the best case scenario, a tablet---usually held by its original shrine or temple, a local museum, or a community center---would be on public display in a well-lit location. Many tablets, on the contrary, are buried within cramped storehouses. Frequently, in rural locations, a sangaku will reside in a small locked shrine building, the sole key-holder of which must be tracked down for access.

\begin{figure}
\begin{center}
\includegraphics[width=5in]{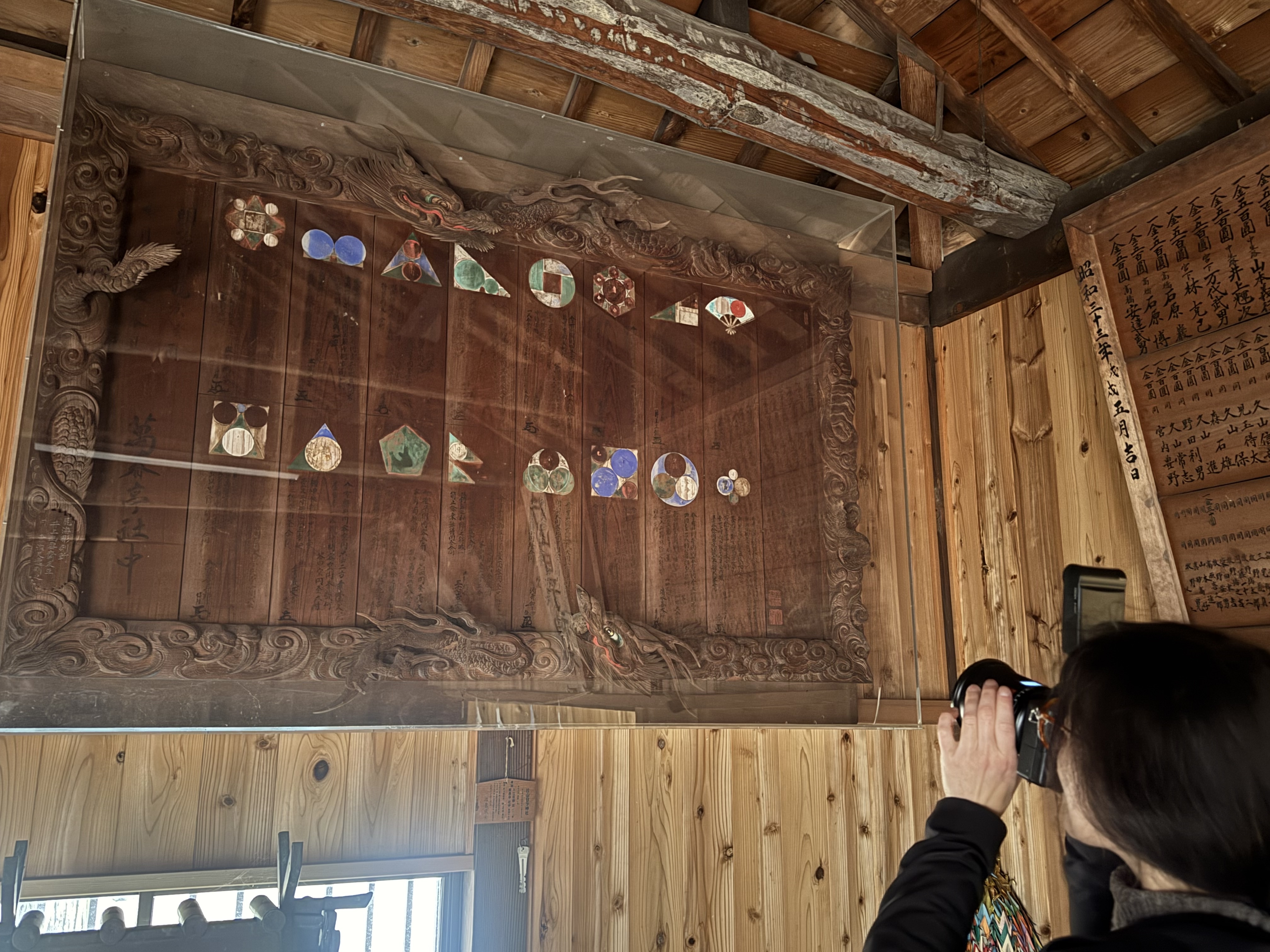}
\caption{Antonia Karaisl photographing a large sangaku in Okayama.}
\end{center}
\end{figure}

Thus, acquiring an image of a single tablet required hours of phone calls, emails, and logistical planning. Japan is well-known for its excellent rail system, but sangaku frequently reside far from train stations. Buses and bicycles filled the gaps. Various obstacles delayed the process. For example, after traveling to Niigata Prefecture in winter, Karaisl and her local contact were unable to enter a tablet-containing shrine after a heavy snow blocked the door. During each successful photoshoot, she discussed the project with the monk, priest, local curator, or village elder responsible for the tablet, establishing permission to publish a composite photographic image online. Many (but not all) sangaku images may be reused under certain conditions (e.g., via a Creative Commons license); this information, of course, is included in each tablet's Archive metadata.

Beyond these caretakers, the phenomenon of sangaku tablets is not well-known in Japan. Inquire about them to a typical Tokyoite and they will most likely think you are asking about mountains
\begin{CJK}{UTF8}{min}
(山岳)
\end{CJK}
 or production amounts
 \begin{CJK}{UTF8}{min}
(産額),
\end{CJK}
two common homophones. I have found myself on the outskirts of conversations in which a shrine attendant thought we were asking to see an Edo-period grain harvest ledger. But within certain circles sangaku are not only familiar; they are regarded as the merely the tip of the wasan iceberg.

\bigskip
\centerline{\rule{2in}{0.4pt}}
\bigskip

On a cool spring evening in 2025---the cherry blossoms in Osaka, just past peak, had begun their slow descent streetward---I sat in a generic classroom in a satellite campus building. The other twenty or so attendees were predominantly older Japanese men, and they arrived with purpose: notebooks and pens at the ready, many were quietly studying handouts. This was a meeting of the Kinki Wasan Research Group.
Named for Japan's Kinki region (comprising seven prefectures in central Honshu), the longstanding 
organization holds a monthly seminar in Osaka; the gatherings are open to all (even a \emph{gaijin} like myself) and draw enthusiasts, amateur scholars, and academics alike. For the next two hours the room dove deeply into a small handful of wasan geometry problems. Slideshows of diagrams and original documents sparked intense discussion. (The longest, most charismatic presentation was given by one Hiroshi Kotera\footnote{Dr. Kotera sadly passed away on July 19, 2025.}.) The topics were rigorous, but laughter was common. The seminar participants---most of whom were there as an act of pure recreation---were having fun. Their keen attention, though, also struck me as a form of reverence. After all, these centuries-old problems were bequeathed to them. 

Other groups like this exist throughout Japan;
for example, there are active wasan research organizations in Nagano, Fukushima, Iwate, and Yamagata Prefectures. Among many valuable functions, these groups frequently maintain detailed records of the sangaku tablets posted in their region, including translations of the problems into modern Japanese. This is nontrivial: like other scholarly documents of the time, traditional sangaku were written in Classical Chinese (\emph{kanbun} or
\begin{CJK}{UTF8}{min}
漢文),
\end{CJK}
making them difficult to parse without specialized knowledge in both language and mathematics. (It was an early hope to include these translations, along with their English equivalents, on the Archive; unfortunately, the resources necessary were not available.) Still, many sangaku problems are available in English translation, some also with their historical solutions (see \cite{fukped,fukrig,fukrot2,hosking2,unger30}). These are an appealing entry point into a broader study of wasan, which has a fascinating history of individuals, obstacles, and accomplishments.

If traditional Japanese mathematics were a garden, it would have high walls and Chinese roots. Early on, Japan's big neighbor to the west gave them physical tools---the abacus and a system of calculating rods, known respectively in Japan as \emph{soroban}
\begin{CJK}{UTF8}{min}
(算盤)
\end{CJK}
and \emph{sangi}
\begin{CJK}{UTF8}{min}
(算木)
\end{CJK}
---and techniques for numerically finding the roots of polynomials and calculating $\pi$, as well as a large collection of geometrical theorems. Not long into the Edo period, however, shogun Tokugawa Iemitsu's isolationist policy of \emph{sakoku}
\begin{CJK}{UTF8}{min}
(鎖国,
\end{CJK}
literally ``locked country") slowed cultural and intellectual exchange with the outside world to a trickle. 
Thus, wasan developed nearly completely devoid of European influence; the importation of a Dutch trig table in the mid-17th century was a rare (and minimally impactful) exception \cite{unger3}. Otherwise, the soil, sun, and rain that nourished wasan's garden came from within. Lacking Western tools like coordinate geometry and differential calculus, Japanese mathematicians made real progress on hard problems. Among other feats, they found and used infinite series expansions for many functions; developed a form of integral calculus; and generated a robust set of techniques for solving systems of equations. (See \cite{smith1914,MR130808,horiuchi1,MR1893990,unger-wedge,clark3}.)
Many results typically attributed to Europeans were independently discovered by Edo period mathematicians, and examining their approaches to problems---in contrast with modern methods---can be deeply insightful.

\bigskip
\centerline{\rule{2in}{0.4pt}}
\bigskip

The Sangaku Archive is a cultural preservation project that also has great potential as a resource for both scholarship and pedagogy. In addition to its database of images and metadata, the Archive hosts a map with over a thousand pins, each of which marks a tablet's location.
These features will aid historians in exploring the social networks of Japanese mathematicians (e.g., master-disciple relationships, routes of knowledge transmission, and feuds), for which the sangaku tradition provided a means of connection. (See \cite{MR4331774,wong,unger3} for pre-Archive work of this nature.)

For teaching geometry and math history, sangaku are an excellent hook: these unique, multilayered objects give math problems a physical presence in time and space, and connect them to the fascinating broader story of wasan.
I have led dozens of wasan workshops for secondary teachers and students, and taught college courses for math-majors and non-majors alike; in my experience, weaving mathematics together with history, art, religion, and Japanese culture makes that math significantly more appealing. Add to this the fact that, for geometry enthusiasts, many sangaku problems are simply fun to work on, and possibly quite different from what they have done in secondary school.

As mentioned earlier, the Sangaku Archive itself does not include translated problems for each tablet. It will, however, have a companion library of Geogebra modules for selected problems. Geogebra---a free, versatile, online math tool---allows users to create geometric figures with compass-and-straightedge (and other) constructions; labels, text, sliders, colors, and animation can also be incorporated. The interactive modules paired with the Archive will
include beginner-level problems with scaffolded hints; harder problems with discussions of solutions; and useful theorems (e.g., the Crossed Chords Theorem and the Descartes Circle Theorem), many of which were known (in some form) to Japanese mathematicians. When appropriate, modules will include links back to related Archive entries, as well as relevant historical notes and bibliographic references.

\begin{figure}
\begin{center}
\fbox{\includegraphics[width=5in]{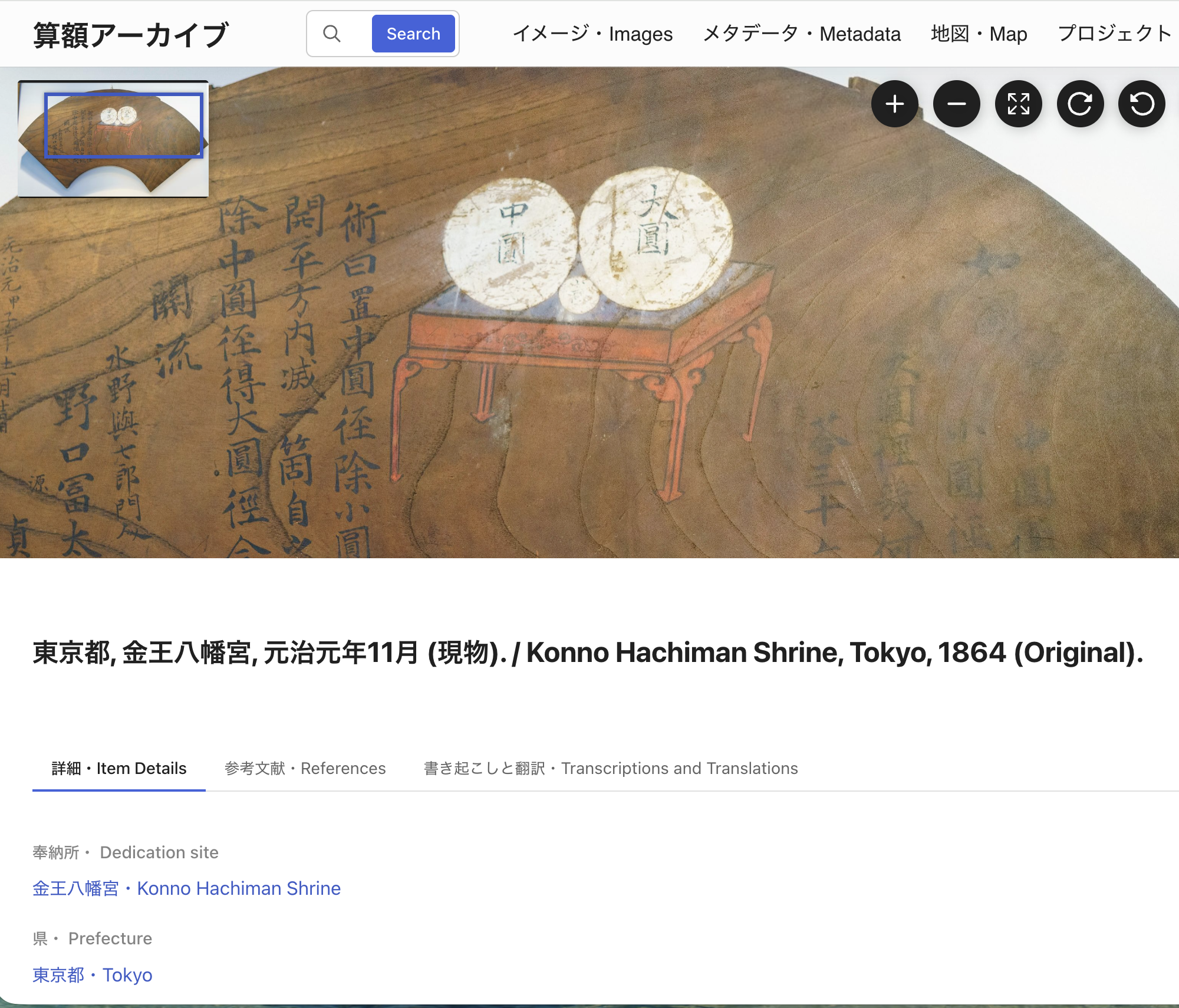}}

\vspace{.15in}

\fbox{\includegraphics[width=5in]{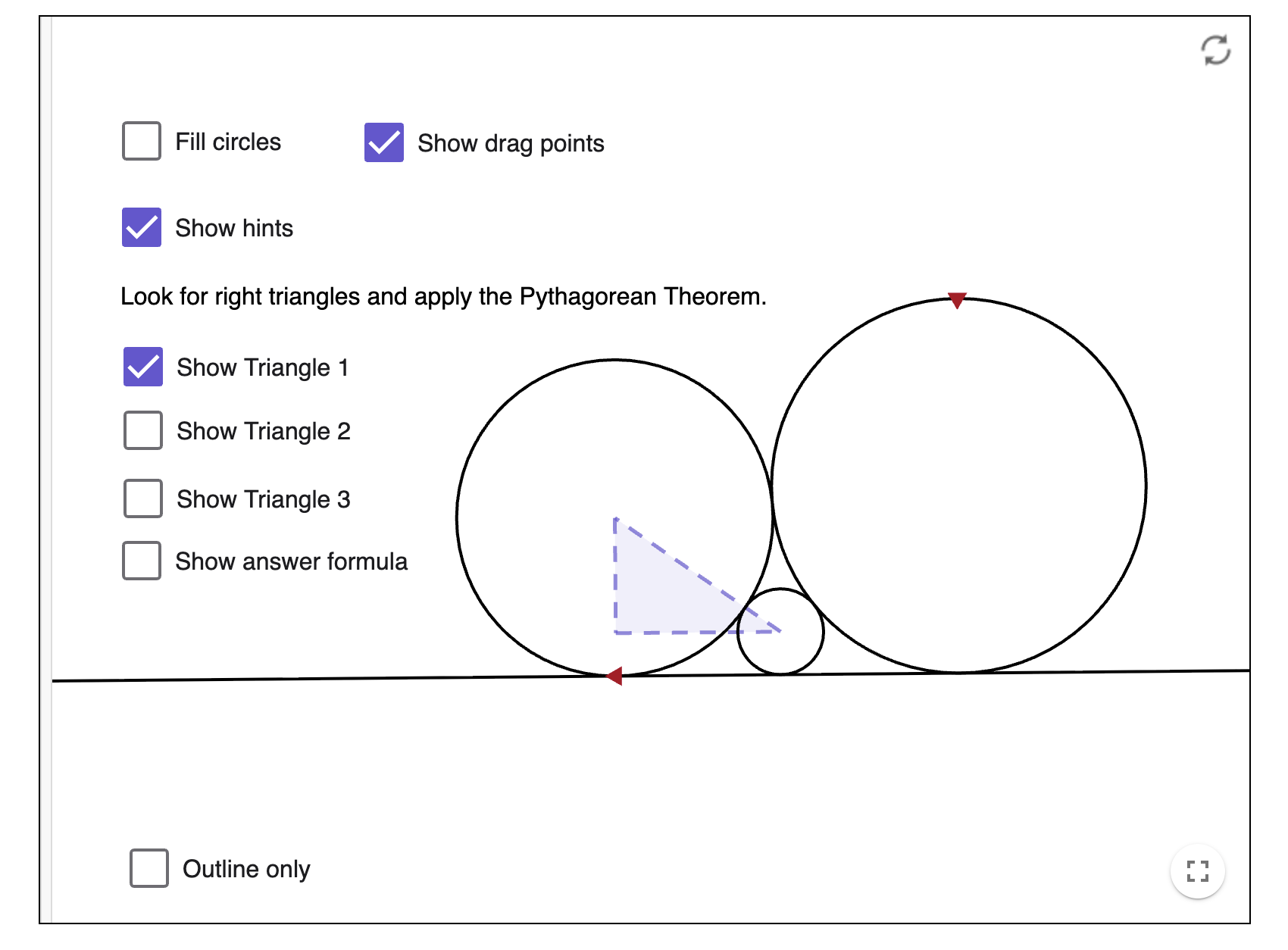}}
\caption{Screenshots of an Archive entry and its corresponding Geogebra module.}
\end{center}
\end{figure}

Like the sangaku tradition, Geogebra has features that foster social networking. Modules and activities can be saved, shared, and compiled, and you can ``follow" other contributors. (I am currently following a Spanish mathematician who has made many gorgeous modules, including a few of sangaku problems.) Our hope is that Geogebra will become an online gathering place for sangaku educators, students, and enthusiasts. While the modules currently paired with the Archive represent a tiny fraction of all sangaku problems (and even tinier of wasan more generally), such a community could build and add their own Geogebra creations. This would expand the set of explorable problems online
(while removing barriers for folks who don't have access to the books and papers that hold English translations), and also generate interest in the Archive, giving users historical context for the mathematics and encouraging further investigation. The seeds of wasan's garden could thus be spread far and wide.

\bigskip
\centerline{\rule{2in}{0.4pt}}
\bigskip

In 1872, the newly-created Ministry of Education issued an edict forbidding the teaching of wasan. Less than a decade earlier, with the fall of the Tokugawa shogunate and forcible reopening of Japan, Western methods of math and science began flooding into the country. Reforms were deemed necessary \cite{ravina}.

Compliance with the edict was neither immediate nor universal. Rural schools, away from the centers of power, were slow to adopt the Western math (also known as \emph{yosan} or
\begin{CJK}{UTF8}{min}
洋算)
\end{CJK}
that would eventually supplant wasan entirely. 
As these tectonic shifts were profoundly changing Japan, the tradition of dedicating sangaku tablets persisted. According to the Archive, more than a hundred tablets
were hung during the Meiji period (1868--1912) alone. 
\begin{figure}
\begin{center}
\includegraphics[width=5in]{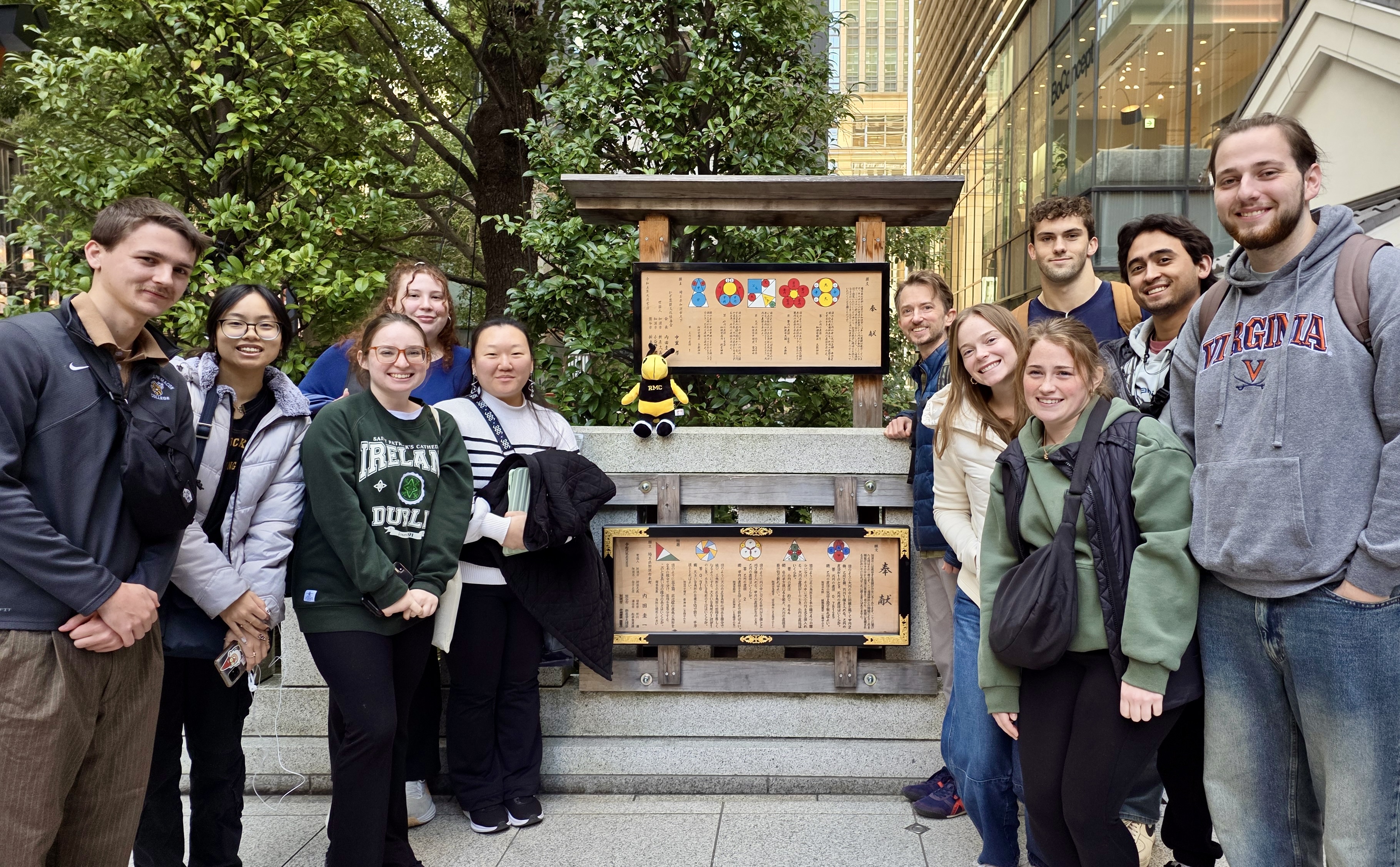}
\caption{The author and his students visiting two new sangaku in Tokyo.}
\end{center}
\end{figure}

Since then, the pace of sangaku dedications has slowed, but it has not stopped. When I took students to Japan in January 2026, we visited three very recent additions to the corpus. One of them is a bilingual tablet of new problems dedicated in 2018 by New Zealander Rosalie Hosking\footnote{Dr. Hosking's Ph.D. thesis on sangaku was referenced earlier.}; her sangaku is displayed at Kyoto's Kitano Tenmangu shrine under the same roof as one from the Meiji era. (Given the relative scarcity of accessible tablets in major cities, such ``two-for-one" sites are the
\emph{pi\`eces de r\'esistance} of my travel course.)
 
 There are also two new tablets at Fukutoku shrine in Tokyo, just outside the Nihombashi subway station. Both were dedicated within the last six years by the Kazo Sangaku Cultural Preservation Society and display a curated selection of problems---translated into modern Japanese---from Edo period tablets hung in and around Kazo City in Saitama Prefecture. The Archive has documented about fifty
replica-style tablets like this, though most are strict reproductions of individual sangaku, unlike the composites at Fukutoku.

The creation of a new tablet---replica or otherwise---is not a minor undertaking. After the mathematics is finalized, artisans hand-paint the text and figures on polished Japanese cypress using historically matched colors; the frame requires precise carpentry and custom hardware. Parts and labor for a single tablet can cost upwards of \$10K. We learned all of this from Takashi Nakazato, a junior high school teacher, Preservation Society member, and driving force behind one of the Fukutoku sangaku. He told us the Society aims ``to pass down this valuable local culture to future generations, creating a homeland where our descendants 100 or 150 years from now can feel proud of the town they were born and raised in, just as we do today."


This veneration of the sangaku tradition impressed strongly upon my students---how many of us have seen mathematics cast in such a light before? Our students may not feel such ancestral connections to geometry, but we can guide them to tap into that energy. Besides, mathematics transcends cultural boundaries. These sangaku problems have been bequeathed to us, too, and we now have the modern digital tools to help us learn from them.

\bigskip

\noindent Please note: the projects below are WORKS IN PROGRESS.

\bigskip

\noindent Visit the Sangaku Archive: \url{https://sangaku-archive.org/}

\bigskip

\noindent Visit the companion Geogebra library: \url{https://www.geogebra.org/m/vp4rbsjk}


\newpage
\bibliographystyle{plain}
\bibliography{references}

\end{document}